\documentstyle{article}
\newtheorem{theorem}{Theorem}[section]
\newtheorem{pro}{Proposition}[section]

\newtheorem{remark}{Remark}[section]

\newcommand{\proof}[1]{\noindent{\it\bf Proof:#1\ }}
\newcommand{\QED}{\hfill$\Box$\medskip}

\begin{document}

\title{  $C^{m_0}$-Smoothness of Evaluation Maps}
\author{
Gang Liu }
\date{November 25,  2013}
\maketitle

In this note, we give  a complete proof of the  claim in [L2] that  total $l$-fold evaluation map $E: \Sigma^l\times {\cal M}\rightarrow M^l$ given by $E(x, f)=(f(x_1), \cdots, f(x_l))$  is of class  $C^{m_0}$.  Here $\Sigma$ and $M$  are two   $C^{\infty}$-smooth Riemannian  manifolds with $\dim (\Sigma)=n, $
${\cal M}={\cal M}_{k, p}$ is the space of $L_k^p$-maps from $\Sigma$ to $M$ and $m_0$ is a positive integer such that $m_0+\gamma=k-\frac{n}{p}$ with  $0<\gamma<1$. Note that $m_0$  is the Sobolev differentiability of any generic elements in ${\cal M}.$

The proof is written only using elementary calculus on Banach spaces.

\begin{theorem}
The total evaluation map $ E: \Sigma^l\times {\cal M}\rightarrow M^l$ is of class $C^{m_0}$.

\end{theorem}

\proof

It is easy to see that $E$ is linear,  and hence, of $C^{\infty}$ in ${\cal M}$-direction.  It is of $C^{m_0}$ in $\Sigma$-direction by Sobolev embedding. The question is about those
mixed partial derivatives as well as the continuity of all derivatives.  

To this end,  we make a few reductions. Clearly, this  can be reduced to the case that $l=1,$  and it can be reduced further  first to  the case that  $M={\bf R}^k$ by using an  embedding of $M$ into ${\bf R}^k$, then to the case that $M={\bf R}^1$. 
 
Since the computations for the partial derivatives are local in $\Sigma$, by multiplying a fixed cut-off function on $\Sigma$ supported near the point that we are interested,  we may assume that $\Sigma$ is either ${\bf R}^n$ or ${\bf T}^n$, the  $n$-tours.  In this setting, $E$ becomes $E:
\Sigma\times L_k^p(\Sigma, {\bf R}^1)\rightarrow {\bf R}^1.$ We will denote $L_k^p(\Sigma, {\bf R}^1)$ by $L_k^p$.

To compute these partial derivatives, let ${\cal D}$ be the space of all smooth function on $\Sigma$ (with compact support if $\Sigma$ being ${\bf R}^n$) with "$C^{\infty}$-topology" in the sense of distribution theory, and ${\cal D}'$ is the collection of continuous linear functionals.    
Consider  the collection of those elements of ${\cal D}'$ that can be extended to  continuous linear functionals on  $L_k^p$.
We denote it by $(L_k^p)'$  with operator norm.  Note that  $(L_k^p)'$ is usually denoted by $L_{-k}^q$ with $1/p+1/q=1.$ Its elements  have more concrete expressions.  But we will only consider $(L_k^p)'$ as an abstract dual.

  Let $\delta:\Sigma\rightarrow {\cal D}'$ defined by $\delta(x)=\delta_x\in {\cal D}'$ where $\delta_x$ is the Dirac delta function at $x\in \Sigma.$  By our assumption, $\delta(x)$ is in $(L_k^p)' $.  Therefore, we have $\delta:\Sigma\rightarrow {(L_k^p)}'$, and $E$ is the composition of $\delta\times Id_{L_k^p}:\Sigma\times {L_k^p} \rightarrow   {(L_k^p)}' \times  {L_k^p}$ with the paring of  ${(L_k^p)}' $ and $  {L_k^p}.$ 

Now we list  the following three elementary facts proved, for instance,  in Lang's book "Real Analysis":

\vspace{2mm}
\noindent (I) Any  paring, as a bilinear continuous map $<-,- >:E_1\times E_2\rightarrow E_3$ between Banach spaces satisfying the condition that 
\vspace{2mm}
\noindent
$\|<e_1, e_2>\|\leq \|e_1\|\cdot \|e_2\|$, is of class $C^{\infty}.$ 

\vspace{2mm}
\noindent(II) A map $F=F_1\oplus F_2:E\rightarrow E_1\oplus E_2$ between Banach  spaces is $C^r$-smooth if and only if each $F_i, i=1, 2,$ is. 

\vspace{2mm}
\noindent(III) The projection $p_i: E_1\oplus E_2\rightarrow E_i, i=1, 2,$ is linear, and hence $C^{\infty}$-smooth.


\vspace{2mm}
\noindent
Note that in our case, for the paring $<-, ->: (L_k^p)'\times L_k^p \rightarrow {\bf R}^1$, we have 
$$|<\phi, \xi>|=|\phi( \frac {\xi}{\|\xi\|_{k, p}})|\cdot \|\xi\|_{k, p}\leq \sup_{\|\eta\|_{k, p}\leq 1}|\phi( \eta)|\cdot \|\xi\|_{k, p}$$
$$=\|\phi \|_{(L_k^p)'}\cdot \|\xi\|_{L_k^p}.$$

Using the above  three facts,  we only need to show that $\delta$ is of $C^{m_0}$. To this end, we observe that for each $x \in \Sigma,$ $\delta_x$ extends to $L_{k-m_0}^p$ since by our assumption  $L_{k-m_0}^p$ is in $C^{\gamma}$ with $0<\gamma<1.$  In other words,  the map $\delta_x$  is lifted as $\delta_x:\Sigma\rightarrow (L_{k-m_0}^p)'\subset (L_{k}^p)'.$
 The following fact will be used repeatedly:  for any $\xi$ in $(L_{k-m_0+l}^p)'$ with $l\leq m_0, $ $\|\xi\|_{(L_{k}^p)'}\leq \|\xi\|_{(L_{k-m_0+l}^p)'},$   which follows from the dual version of the inequality.

 The result we are looking for follows from this observation.  Roughly speaking, each time we take a partial derivative to $\delta$, we move it from the dual of $L_{k-i}^p$ to the dual  of
$L_{k-i+1}^p$ starting with $i=m_0$.

More precisely, we show  this inductively by the following four steps:

\medskip
\noindent
${\bullet}$ Step I:
$\delta:\Sigma\rightarrow   (L_{k-m_0}^p)'$ is continuous with respect to the operator norm on $(L_{k-m_0}^p)'$.

\proof

$$||\delta (x)-\delta (y)||_{(L_{k-m_0}^p)'}=sup_{||\xi||_{k-m_0, p}\leq 1}|| (\delta (x)-\delta (y))(\xi)||$$$$=  sup_{||\xi||_{k-m_0, p}\leq 1}|| \xi(x)-\xi(y)||$$ $$\leq sup_{||\xi||_{k-m_0, p}\leq 1}|| \xi||_{C^{0, \gamma}}||x-y||^{\gamma}$$$$\leq C_0\cdot sup_{||\xi||_{k-m_0, p}\leq 1}|| \xi||_{k-m_0,p}||x-y||^{\gamma}$$$$=C_0\cdot ||x-y||^{\gamma}$$ for some constant $C_0$.

Since $||\delta (x)-\delta (y)||_{(L_{k}^p)'}\leq ||\delta (x)-\delta (y)||_{(L_{k-m_0}^p)'}$, this also proves that $\delta$ is continuous.
   \QED

 \noindent
 ${\bullet}$ Step II: (A )The value of the partial derivative of $\delta$ at $x\in \Sigma,$ $({\partial_j\delta})(x)$ is equal to the distribution derivative of $\delta_x$, ${\partial_j(\delta_x)}.$ (B)${\partial_j(\delta_x)}\in (L_{k-m_0+1}^p)'\subset (L_{k}^p)'$.  Therefore, ${\partial_j\delta}:\Sigma\rightarrow   (L_{k-m_0+1}^p)'$ defined by $({\partial_j\delta})(x)={\partial_j(\delta_x)}.$

 \proof

 Since for any $\xi\in L_{k-m_0+1}^p$, $$||{\partial_j(\delta_x)}(\xi)||=||\delta_x({\partial_j(\xi)})||$$
 $$=||{\partial_j(\xi)}(x)||$$$$\leq
 ||{\partial_j(\xi)}||_{C^0}\leq  C_1\cdot ||{\partial_j(\xi)}||_{k-m_0, p}$$
 $$= C_1\cdot ||\xi||_{k-m_0+1, p}.$$ This shows that (B) is true.

 To prove (A), we compute

 $$||\frac {\delta (x+he_j)-\delta (x)}{h}-{\partial_j(\delta_x)}||_{(L_{k-m_0+1}^p)'}$$
 $$=sup_{||\xi||_{k-m_0+1, p}\leq 1}|| \frac {\delta (x+he_j)-\delta (x)}{h}(\xi)-{\partial}_j(\delta_x (\xi))||$$$$=
  sup_{||\xi||_{k-m_0+1, p}\leq 1}||\frac {\xi (x+he_j)-\xi (x)}{h}-\delta_x({\partial_j\xi})||
  $$$$=sup_{||\xi||_{k-m_0+1, p}\leq 1}||\frac {\xi (x+he_j)-\xi (x)}{h}-{\partial_j\xi}(x)||$$$$
  =sup_{||\xi||_{k-m_0+1, p}\leq 1}||{\partial_j\xi}(x+te_j)-{\partial_j\xi}(x)||
  $$$$\leq sup_{||\xi||_{k-m_0+1, p}\leq 1}||{\partial_j\xi}||_{C^{0, \gamma}}||te_j||^{\gamma}$$
  $$\leq   sup_{||\xi||_{k-m_0+1, p}\leq 1}||{\xi}||_{C^{1, \gamma}}||te_j||^{\gamma}
  $$$$\leq   C_1\cdot sup_{||\xi||_{k-m_0+1, p}\leq 1}||{\xi}||_{k-m_0+1, p}||te_j||^{\gamma} $$$$
  \leq   C_1\cdot ||h||^{\gamma}.$$

  \noindent
  Here $0<t<h$.

 Therefore,   $$||\frac {\delta (x+he_j)-\delta (x)}{h}-{\partial_j(\delta_x)}||_{(L_{k}^p)'}$$
 $$\leq ||\frac {\delta (x+he_j)-\delta (x)}{h}-{\partial_j(\delta_x)}||_{(L_{k-m_0+1}^p)'}$$  $$\leq   C_1\cdot ||h||^{\gamma}\rightarrow 0$$ as $h\rightarrow 0.$ This proves that ${\partial}_j\delta$ exists.

  \QED

  \noindent
  ${\bullet}$ Step III:  Assume that ${\partial^{\alpha} \delta}:\Sigma\rightarrow   (L_{k-m_0+l}^p)'\subset (L_{k}^p)'$ for multi-indices $\alpha=(\alpha_1, \cdots,\alpha_n)$ with $|\alpha|=l\leq m_0, $ and ${\partial^{\alpha} \delta}(x)={\partial^{\alpha} (\delta_x)}$.
 Then  ${\partial^{\alpha} \delta}:\Sigma\rightarrow   (L_{k-m_0+l}^p)'\subset (L_{k}^p)'$ is continuous.

 \proof

  $$||{\partial^{\alpha}\delta} (x)-{\partial^{\alpha}\delta} (y)||_{(L_{k-m_0+l}^p)'}$$$$=sup_{||\xi||_{k-m_0+l, p}\leq 1}||
  ({\partial^{\alpha}\delta} (x)-{\partial^{\alpha}\delta} (y))(\xi)||$$$$=  sup_{||\xi||_{k-m_0+l, p}\leq 1}
  || {\partial^{\alpha}(\delta_x)}(\xi) -{\partial^{\alpha}(\delta_y)}(\xi)||$$
  $$= sup_{||\xi||_{k-m_0+l, p}\leq 1}
  || {\partial^{\alpha}(\xi)}(x) -{\partial^{\alpha}(\xi)}(y)||
  $$$$\leq  sup_{||\xi||_{k-m_0+l, p}\leq 1}
  || {\partial^{\alpha}(\xi)}||_{C^{0, \gamma}}||x-y||^{\gamma}
  $$$$\leq  sup_{||\xi||_{k-m_0+l, p}\leq 1}
  || {\partial^{\alpha}(\xi)}||_{k-m_0,p}||x-y||^{\gamma}$$$$\leq C_l\cdot sup_{||\xi||_{k-m_0+l, p}\leq 1} || \xi)||_{k-m_0+l,p}||x-y|^{\gamma}
  $$$$ =C_l\cdot ||x-y||^{\gamma}
  $$

 \QED

\noindent
${\bullet}$ Step IV:  For multi-indices $\alpha=(\alpha_1, \cdots,\alpha_n)$ with $|\alpha|=l\leq m_0-1, $ assume that ${\partial^{\alpha} \delta}:\Sigma\rightarrow ( L_{k-m_0+l}^p)' \subset (L_{k}^p)'$  is continuous. Then  (A) $({\partial}_j{\partial^{\alpha} \delta})(x)={\partial}_j({\partial^{\alpha}} (\delta_x))$; and (B)
 ${\partial}_j({\partial^{\alpha}} (\delta_x))$ is in $(L_{k-m_0+l+1}^p)'\subset (L_{k}^p)'.$

 \proof

 Since for any $\xi\in L_{k-m_0+l+1}^p$, 
 $$||{\partial}_j({\partial^{\alpha}(\delta_x)})(\xi)||$$
 $$=||\delta_x({\partial}_j{\partial}^{\alpha}(\xi))||=||({\partial}_j{\partial}^{\alpha}(\xi))(x)||$$
 $$\leq
 ||{\partial}_j{\partial}^{\alpha}(\xi)||_{C^0}$$
 $$ \leq  C_{l+1}\cdot ||{\partial}_j{\partial}^{\alpha}(\xi)||_{k-m_0, p}$$
 $$= C_{l+1}\cdot ||\xi||_{k-m_0+l+1, p}.$$ This shows that (B) is true.

 To prove (A), we compute

 $$||\frac {{\partial}^{\alpha}\delta (x+he_j)-{\partial^{\alpha}\delta} (x)}{h}-{\partial}_j {\partial}^{\alpha}(\delta_x)||_{(L_{k-m_0+l+1}^p)'}$$$$=
 sup_{||\xi||_{k-m_0+l+1, p}\leq 1}|| \frac {{\partial}^{\alpha}\delta (x+he_j)-{\partial^{\alpha}\delta} (x)}{h}(\xi)-{\partial}_j{\partial^{\alpha}(\delta_x)}(\xi)||$$     $$=
  sup_{||\xi||_{k-m_0+l+1, p}\leq 1}||\frac {{\partial^{\alpha}\xi} (x+he_j)-{\partial^{\alpha}\xi} (x)}{h}-{\partial}_j{\partial^{\alpha}\xi}(x)||$$     $$ =sup_{||\xi||_{k-m_0+l+1, p}\leq 1}||{\partial}_j{\partial^{\alpha}\xi}(x+te_j)-{\partial}_j{\partial^{\alpha}\xi}(x)||$$ 
      
$$\leq 
sup_{||\xi||_{k-m_0+l+1, p}\leq 1}||{\partial}_j{\partial}^{\alpha}\xi||_{C^{0, \gamma}}||te_j||^{\gamma}$$     
 $$\leq  C_{l+1}\cdot sup_{||\xi||_{k-m_0+l+1, p}\leq 1}||{\partial}_j{\partial^{\alpha}}{\xi}||_{{k-m_0, p}}||te_j||^{\gamma}  $$     
  $$ =   C_{l+1}\cdot sup_{||\xi||_{k-m_0+l+1, p}\leq 1}||{\xi}||_{k-m_0+l+1, p}||te_j||^{\gamma}
  $$     
  $$\leq   C_{l+1}\cdot ||h||^{\gamma}.$$

 \noindent
 Here $0<t<h$.

  \QED

\begin{remark}
In the above computation, we only prove that all partial derivatives of $\delta:\Sigma \rightarrow (L_k^p)'$ exist and are continuous  up to degree $m_0$. Since the domain $\Sigma$  is of finite dimensional, this  is equivalent to  that $\delta$  is of class $C^{m_0}$ in the usual sense of the differential calculus in Banach
 space setting (see Lang's "Real Analysis" for the proof of this equivalency).  In particular, the proof here has nothing to do with the $sc$-smoothness in the usual polyfold theory.
\end{remark}

\medskip
\vspace{2mm}
\noindent $\bullet$ Note: Proposition 3.1 in [L2] is a corollary of the above theorem, which we state now.

 \begin{pro}
Let $G$ be a Lie subgroup of the group of differomorphisms of $\Sigma$. Fix an $x$ in $\Sigma$. Let $\Phi_x:G\times L_k^p(\Sigma, M)\rightarrow L_k^p(\Sigma, M)\rightarrow M$  be the composition  $ev_x\circ \Psi.$ Here $ev_x:L_k^p(\Sigma, M)\rightarrow M$ is the evaluation map at $x$ and $ \Psi:G\times L_k^p(\Sigma, M)\rightarrow L_k^p(\Sigma, M)$ is the total action map of $G$ acting  on $L_k^p(\Sigma, M)$ as reparametrization group of $\Sigma$.
Then $\Phi_x$ is of class $C^{m_0}.$

\end{pro}

\proof

For the completeness, we include the argument in [L2] that reduces this proposition to the above theorem.

For any $g\in G$ and $\xi \in L_k^p(\Sigma, M)$, we have $$\Phi_x(g, \xi )=ev_x\circ \Psi (g, \xi)=ev_x(\xi\circ g)=\xi(g(x))$$ 
$$=E(g(x), \xi)=E(\phi_x(g), \xi).$$ Here $\phi_x:G\rightarrow \Sigma$ is the orbit map of $x \in \Sigma $ given by $\phi_x(g)=g(x)$ which is $C^{\infty}$-smooth by our assumption that $G$ acts on $\Sigma$ smoothly. Therefore, $\Phi_x=E\circ (\phi_x, Id).$ Here $Id$ is the identity map on $L_k^p(\Sigma, M)$, and $(\phi_x, Id): G\times L_k^p(\Sigma, M)\rightarrow  \Sigma \times L_k^p(\Sigma, M)$ is of class $C^{\infty}$.

\QED

As for the smoothness of $E$, the proof in [L2]  only establishes  the trivial fact that $E$ is of class  $C^{\infty}$ along ${\cal M}$-direction and of class $C^{m_0}$  along $\Sigma$-direction.  Even the continuity of the first derivative is not proved in [L2].   The proof  the Theorem 0.1 above is taken from [L1]. It is possible to give a more direct  proof   for  the $C^{r}$-smoothness of $E$ at least for small values of $r$ starting with the continuity of the first derivative.  However, the computation below shows that similar considerations as above proof has to be used. In the following we  carry out this computation for $C^1$-smoothness of $E$. It also gives another way to reduce the proof of the  Theorem 0.1 to the above statement of the $C^{m_0}$-smoothness of the  ${\delta}$-function.

\medskip
\vspace{2mm}
\noindent $\bullet$ Continuity of $E=E(k-m_0):\Sigma \times L_{k-m_0}^p\rightarrow {\bf R^1}.$ 

 By Sobolev embedding, for any $g$ in $L_{k-m_0}^p$, we have $\|g\|_{C^{0, \gamma}}\leq C\|g\|_{k-m_0,p}$ for a positive constant $C$.

$$  |{E(x+v, f+\xi)- E(x, f)}|= |(f(x+v)-f(x))+\xi(x+v)|$$ 
$$\leq \Sigma_{i=0}^{n-1}|(f(x+v^{i+1})-f(x+v^i))|+|\xi(x+v)|$$ 
$$\leq 
\Sigma_{i=0}^{n-1}\|f\|_{C^{0, \gamma}}|v_{i+1}|^{\gamma}+\|\xi\|_{C^{0}}$$ $$\leq C_1\cdot (\|f\|_{k-m_0, p}\|v\|^{\gamma}+\|\xi\|_{k-m_0, p}).$$

Here $v=(v_1, v_2, \cdots, v_n)$ and $v^i=(v_1, v_2, \cdots, v_i, 0, \cdots 0).$ This proves the continuity directly for $E$ extended to 
$\Sigma \times L_{k-m_0}^p$. It will be the starting point for the induction below.

\medskip
\vspace{2mm}
\noindent $\bullet$ $$DE_{x, f}(v, \xi)=Df_x(v)+\xi(x)=\Sigma_{i=1}^n ({\partial}_if(x)\cdot v_i)+\xi(x).$$ Here $v=(v_1, \cdots, v_n)$ is a tangent  vector in $T_x\Sigma\simeq {\bf R}^n$ and $\xi$ is a tangent vector in $T_fL_k^p\simeq L_k^p$.

\proof

$$  {E(x+tv, f+t\xi)- E(x, f)}$$   $$=(f(x+tv)-f(x))+t\xi(x+tv)=(Df_{x+htv}(v)+\xi(x+tv))t.$$ Here $0<h<1$

$$\| \frac {E(x+tv, f+t\xi)- E(x, f)}{t}-(Df_x(v)+\xi(x))\|$$ 
$$\leq \|(Df_{x+htv}-Df_x)(v)\|+\|\xi(x+tv)-\xi(x)\|$$
$$\leq \|(Df_{x+htv}-Df_x)\|\cdot \|v\|+\|\xi(x+tv)-\xi(x)\|$$ 
$$\leq \|Df\|_{C^{0, \gamma}}|ht|^{\gamma}\|v\|^{1+\gamma}+\|\xi\|_{C^{0, \gamma}}|t|^{\gamma}\|v\|^{\gamma}$$ 
$$\leq  \|Df\|_{k-m_0, p}|ht|^{\gamma}\|v\|^{1+\gamma}+\|\xi\|_{k-m_0, p}|t|^{\gamma}\|v\|^{\gamma}.$$ 
$$\leq  \|f\|_{k, p}|ht|^{\gamma}\|v\|^{1+\gamma}+\|\xi\|_{k, p}|t|^{\gamma}\|v\|^{\gamma}\rightarrow 0$$ as $t\rightarrow 0.$
\QED

Therefore, $DE: \Sigma\times L_k^p\rightarrow L({\bf R}^n\oplus L_k^p, {\bf R})$ is given by $DE(x, f)=Df_x+\delta_x$.
Note that $$L({\bf R}^n\oplus L_k^p, {\bf R})\simeq L({\bf R}^n , {\bf R})\oplus L(L_k^p, {\bf R})=({\bf R^n})'\oplus (L_k^p)'.$$

\noindent $\bullet$  Degree of smoothness of $DE$:

\medskip
\vspace{2mm}
\noindent 

 Note that the second term of $DE$ is exactly the map $\delta:\Sigma\times L_k^p\rightarrow  (L_k^p)'\rightarrow L({\bf R}^n\oplus L_k^p, {\bf R})$ before. It is of class $m_0$.

Denote the first term of $DE$ by $E^1:\Sigma\times L_k^p\rightarrow ({\bf R^n})'=(({\bf R}^1)')^n$. Then $E^1$ can  be identified 
with $E^1\simeq (E^1_1, \cdots, E^1_n)$ with $E^1_j:\Sigma\times L_k^p\rightarrow ({\bf R}^1)'= {\bf R}$ defined by $E^1_j(x, f)
={\partial }_jf (x)$. Clearly $E^1_j=E(k-1)\circ (I_{\Sigma} \times P_j(k))$. Here $P_j(k):L_k^p\rightarrow L_{k-1}^p$ is the bounded linear  map given by
$P_j(k)(f)={\partial }_jf$, which is smooth, and $E(k-1):{\Sigma} \times L_{k-1}^p\rightarrow {\bf R}$ is just the total evaluation map $E=E(k)$
 with a shifting from level $k$ to $k-1$.    We already proved that $E(k-m_0)$  is continuous. We are in the position to apply induction to conclude that $E(k-1)$ is of class $m_0-1$ provided that required smoothness for the delta function is already established.   This implies that the  first  term of $DE$ is of  class $m_0-1$,  and $E$ is of class $m_0$.
 In particular, the above argument shows that  for the  proof of  $C^1$-smoothness of $E$, only the continuity of $\delta$ is needed.

\medskip
\vspace{2mm}
\noindent $\bullet$ Acknowledgement:  In [L2], the author stated that:

\vspace{2mm}
\noindent  
(A)  Proposition 3.1 in [L2]  is weaker than the following  statements
that (i) the action map $\Psi:G\times {\cal M}_{k, p}\rightarrow {\cal M}_{k-1, p}$ is of class $C^1$; (ii) inductively 
 $\Psi:G\times {\cal M}_{k, p}\rightarrow {\cal M}_{k-l, p}$ is of class $C^l$ with $l\leq m_0.$

\vspace{2mm}
\noindent 
(B) Above (i) and  (ii) follows from the considerations in the theory of $sc$-smoothness.

\medskip 
\vspace{2mm}
\noindent
I am grateful to  McDuff for pointing out that  there is a difference  between  the usual smoothness and  the $sc$-smoothness in Polyfold theory,
 and the statement (B) needs to be clarified.

\medskip 
\vspace{2mm}
\noindent
 Indeed,  the first derivative appeared in (A) (i) above is just the ordinary derivative even we use $sc$-type of computation. However,  the continuity of the  first derivative in the $sc$-smoothness is measured  in  a weaker topology on $L (E_1, E_2)$ (called strong topology
 in operator theory) rather than in norm topology.

\medskip
\vspace{2mm}
\noindent 

The  author's intention  for (A) is to give another proof for  Proposition 3.1 in [L2]. Since our method  for regularizing  the moduli spaces of $J$-holomorphic curves
 in  [L2] does not use $sc$-smoothness, we will  not discuss the statement (A) further here.
 In [L3],  we will  prove that (A) above is true  in the
 sense of  usual calculus on  Banach Manifolds in the case that $(k-l)-\frac {n}{p}>0$. In other words, in the above situation, the $sc$-smoothness in the standard polyfold  theory is not the optimal result for the purpose here despite of the fact that $sc$-smoothness using the weaker topology is the right choice for 
 various other reasons  in polyfold theory.

 \medskip
\vspace{2mm}
\noindent $\bullet$ A question:
\vspace{2mm}
\noindent 

 The question  still remains that whether  or not 
 Proposition 3.1 in [L2]  is weaker than the two statements in (A) interpreted  in the sense of $sc$-smoothness.
  I am expecting a positive answer  to this question.

\end{document}